\documentclass[oneside,12pt]{amsart}
\usepackage{amsmath, amsfonts,amsthm,times,graphics}

 \makeatletter
\renewcommand*\subjclass[2][2010]{%
  \def\@subjclass{#2}%
  \@ifundefined{subjclassname@#1}{%
    \ClassWarning{\@classname}{Unknown edition (#1) of Mathematics
      Subject Classification; using '2010'.}%
  }{%
    \@xp\let\@xp\subjclassname\csname subjclassname@#1\endcsname
  }%
}

 \makeatother

\theoremstyle{definition}

 \makeatletter
\renewcommand*\subjclass[2][2010]{%
  \def\@subjclass{#2}%
  \@ifundefined{subjclassname@#1}{%
    \ClassWarning{\@classname}{Unknown edition (#1) of Mathematics
      Subject Classification; using '1991'.}%
  }{%
    \@xp\let\@xp\subjclassname\csname subjclassname@#1\endcsname
  }%
}
 \makeatother

\begin{document}

\title[Wolstenholme's theorem:  
Its Generalizations and Extensions...]{Wolstenholme's theorem: 
Its Generalizations and Extensions in the last hundred and fifty years
(1862--2012)}

\author{Romeo Me\v strovi\' c}

\address{Maritime Faculty, University of Montenegro, Dobrota 36,
 85330 Kotor, Montenegro} \email{romeo@ac.me}

{\renewcommand{\thefootnote}{}
\setcounter{footnote}{0}}

\maketitle

  \begin{abstract} In 1862,
150 years ago, J.  Wolstenholme 
proved that for any prime $p\ge 5$ the numerator of the fraction
  $$
1+\frac 12 +\frac 13+\cdots +\frac{1}{p-1}
  $$
written in reduced form is divisible by $p^2$
and that the numerator of the fraction 
  $$
1+\frac{1}{2^2} +\frac{1}{3^2}+\cdots +\frac{1}{(p-1)^2}
  $$
written in reduced form is divisible by $p$.

The first of the above congruences, the so-called {\it Wolstenholme's 
theorem}, is a fundamental  congruence in Combinatorial 
Number Theory. In this article, consisting of 11 sections, 
we provide a historical survey of Wolstenholme's type congruences, 
related problems and conjectures. Namely, we present and compare several 
generalizations and extensions 
of Wolstenholme's theorem obtained in the last hundred and fifty years. 
In particular, we present  about  80  variations 
and generalizations of  this theorem including congruences 
for Wolstenholme primes. These congruences are 
discussed here by 33 remarks.
   \end{abstract} 

\section{Introduction}

Congruences modulo primes have been widely investigated since the time
of Fermat. {\it Let $p$ be a prime. Then by  Fermat little theorem, for each 
integer $a$ not divisible by $p$}
  $$
 a^{p-1}\equiv 1\pmod{p}.
  $$
Furthermore, by {\it Wilson theorem, for any prime $p$}
  $$
(p-1)!+1\equiv 0\pmod {p}.
  $$
From Wilson theorem it follows immediately that 
$(n-1)!+1$ is divisible by $n$ if and only if $n$ is a 
prime number. 

"In attempting to discover some analogous expression 
which should be divisible by $n^2$, 
whenever $n$ is a prime, but not divisible if $n$ is a
composite number",  in 1819 Charles Babbage \cite{bab} is led 
to the congruence 
  \begin{equation*}
{2p -1\choose p-1}\equiv 1\pmod{p^2}
  \end{equation*}
for primes $p\ge 3$. In 1862 J. Wolstenholme 
proved that the above congruence holds modulo $p^3$ 
for any prime $p\ge 5$. 

As noticed in  \cite{gr}, many great mathematicians of the nineteenth century
considered problems involving binomial coefficients modulo a prime 
power (for instance Babbage, Cauchy, Cayley, Gauss, Hensel,
Hermite, Kummer, Legendre, Lucas, and Stickelberger). They discovered 
a variety of elegant and surprising theorems which are
often easy to prove. For more information on these classical  results, their
extensions,  and new results about this subject, see Dickson \cite{d},
Granville \cite{gr} and Guy \cite{gu}.

Suppose that a prime $p$ and pair of integers  $n\ge m\ge 0$ 
are given. A beautiful {\it theorem of E. Kummer} from year 1852 
(see \cite{k} and  \cite[p. 270]{d}) states that if $p^r$ is the highest power of $p$ 
dividing ${n\choose m}$, then $r$ is equal to the number of carries when
adding $m$ and $n-m$ in base $p$ arithmetic.
 If $n=n_0+n_1p+\cdots +n_sp^s$ and
$m=m_0+m_1p+\cdots +m_sp^s$ are the $p$-{\it adic expansions} of $n$ and $m$
(so that $0\le m_i,n_i\le p-1$ for each $i$),
then by {\it Lucas's theorem} from year 1878 
(\cite{l}; also see \cite[p. 271]{d} and \cite{gr}), 
 $$
{n\choose m}\equiv \prod_{i=0}^{s}{n_i\choose m_i}\pmod{p}.
  $$
This immediately yields
 \begin{equation}\label{con1}
{np \choose mp}\equiv {n \choose m}\pmod{p}
  \end{equation}
since the same products of binomial coefficients are formed on 
the right side of Lucas's theorem in both cases, other than an extra
${0\choose 0}=1$.

{\it Remark} 1. A direct proof of the  congruence (\ref{con1}), 
 based on a polynomial method, is given in  \cite[Solution of Problem A-5, 
p. 173]{pu}. \hfill $\Box$

Notice that the  congruence (\ref{con1}) with $n=2$ and $m=1$ becomes 
   $$
{2p \choose p}\equiv 2\pmod{p},
  $$
whence by the identity ${2p\choose p}=2{2p-1\choose p-1}$
it follows that {\it for any prime $p$}
  \begin{equation}\label{con2}
{2p -1\choose p-1}\equiv 1\pmod{p}.
  \end{equation}
As noticed above, in 1819 Babbage  (\cite{bab}; also 
see \cite[Introduction]{gr} or \cite[page 271]{d})
showed that the congruence (\ref{con2}) holds modulo $p^2$,
that is, {\it for a prime $p\ge 3$ holds}
  \begin{equation}\label{con3}
{2p -1\choose p-1}\equiv 1\pmod{p^2}.
  \end{equation}
The congruence (\ref{con3}) is generalized 
in 1862 by Joseph Wolstenholme  \cite{w} as presented in the next section.
Namely, {\it Wolstenholme's theorem} asserts that   
  \begin{equation*}
{2p -1\choose p-1}\equiv 1\pmod{p^3}
  \end{equation*}
for all primes $p\ge 5$. 
 Wolstenholme's theorem plays a fundamental role  in Combinatorial 
Number Theory. In this article, we provide a historical survey of 
Wolstenholme's type congruences, related problems and conjectures
concerning to the Wolstenholme primes. 
This  article consists of 11 sections in which we present numerous
generalizations and extensions of Wolstenholme's theorem established 
in the last hundred and fifty years. 

The article is organized as follows. In Section 2, we present extensions of Wolstenholme's theorem
up to modulus $p^9$. In the next section, 
many of these congruences are expressed in terms of Bernoulli numbers.
Section 4 is devoted to the Wolstenholme's type harmonic series congruences.
Certain Wolstenholme's type supercongruences are given in the next section.
In Section 6, we present Ljunggren's congruence and   
Jacobsthal-Kazandzidis congruence and their variations modulo higher 
prime powers. In the next section, we give several 
characterizations of Wolstenholmes primes and related conjectures.
Wolstenholme's type theorems for composite moduli 
are established in Section 8. The converse of Wolstenholme's theorem 
is discussed in Section 9. In the next section, we present some recent 
congruences for binomial sums closely related to Wolstenholme's theorem.
Finally, some $q$-analogues of Wolstenholme's type congruences
are given in the last section of this survey article. 
    
The Bibliography of this article
contains  107 references 
consisting of 12 textbooks and monographs, 90  papers, 
3 problems, Sloane's On-Line Encyclopedia of Integer 
Sequences  and one Private correspondence. In this article, 
some  results of these references 
are cited  as generalizations of  
certain Wolstenholme's type congruences, but  without the expositions 
of related congruences. The total number of citations 
given here is 197.

\section{Wolstenholme's theorem and its extensions up to modulo $p^7$}

In 1862, 150 years ago, at the beginning of his celebrated paper 
"{\it On certain properties of prime numbers}" \cite[page 35]{w},
J. Wolstenholme wrote:

"{\it The properties I propose to prove in this article, 
for any prime number $n>3$, are $(1)$ that the numerator of the fraction
  $$
1+\frac 12 +\frac 13+\cdots +\frac{1}{n-1}
  $$
when reduced to its lowest terms is divisible by $n^2$, $(2)$ the 
numerator of the fraction 
  $$
1+\frac{1}{2^2} +\frac{1}{3^2}+\cdots +\frac{1}{(n-1)^2}
  $$
is divisible by $n$, and $(3)$ that the number of combinations of
$2n-1$ things, taken $n-1$ together, diminshed by $1$, is divisible
by $n^3$.  I discovered the last to hold, for several cases, 
in testing numerically a result of certain investigations, and
after some trouble succeded in proving it to hold universally. The method
I employed is somewhat laborious, and I should be glad if some 
of your readers would supply a more direct proof....}"   

More precisely, the first mentioned 
result of J. Wolstenholme \cite{w} 
asserts that  {\it if $p\ge 5$ is a prime, then the numerator of  the fraction 
 $$
1+\frac{1}{2}+ 
\frac{1}{3}+\cdots +\frac{1}{p-1}
$$ 
written in the reduced form is divisible by $p^2$}. For a proof, 
also see \cite[p. 89]{hw}, \cite[p. 116]{apo} and External 
Links listed in Appendix A). 

From this congruence it can be easily  deduced that 
{\it the binomial coefficient 
${2p -1\choose p-1}$ satisfies the congruence 
 \begin{equation}\label{con4}
{2p -1\choose p-1}\equiv 1\pmod{p^3}.
  \end{equation}
for any prime $p\ge 5$} (see e.g., \cite[p. 89]{hw}, \cite[p. 116]{apo} 
and \cite{bau}). 

As usual in the literature, in this note 
the congruence (\ref{con4}) is also called {\it Wolstenholme's 
theorem}. Notice also that from the identity 
${2n\choose n}=2{2n-1\choose n-1}$, $n=1,2,\ldots,$ 
we see that (\ref{con4}) also may be
written as 
$$
{2p \choose p}\equiv 2\pmod{p^3}.
 $$   

The congruence (\ref{con4}) is generalized by J.W.L. Glaisher 
in 1900. Namely, by a special case of 
Glaisher's congruence (\cite[p. 21]{gl1}, \cite[p. 323]{gl2};   
also cf. \cite[Theorem 2]{m}), {\it for any prime $p\ge 5$ we have} 
    \begin{equation}\label{con5}
{2p-1\choose p-1}\equiv 1-2p \sum_{k=1}^{p-1}\frac{1}{k}\pmod{p^4}.
      \end{equation}
In 1995 R.J. McIntosh \cite[p. 385]{m}  established a generalization 
of (\ref{con5}) modulo $p^5$; he showed that {\it for any prime $p\ge 7$}
 \begin{equation}\label{con6}
{2p-1\choose p-1}\equiv 1-p^2 \sum_{k=1}^{p-1}\frac{1}{k^2}\pmod{p^5}.
      \end{equation}
On the other hand, as an immediate consequence of a result by J. Zhao
in 2007 \cite[Theorem 3  with $n=2$ and $r=1$]{z2}, {\it for any prime
$p\ge 7$}
   \begin{equation}\label{con7}
{2p-1\choose p-1}\equiv 1+2p \sum_{k=1}^{p-1}\frac{1}{k}\pmod{p^5}.
      \end{equation}
In 2010 R. Tauraso  \cite[Theorem 2.4]{t} proved that {\it for any prime
$p\ge 7$}
   \begin{equation}\label{con8}  
{2p-1\choose p-1}\equiv 1 +2p\sum_{k=1}^{p-1}\frac{1}{k}
+\frac{2p^3}{3}\sum_{k=1}^{p-1}\frac{1}{k^3}\pmod{p^6} 
    \end{equation}
{\it which also can be written as \cite[Corollary 1.4]{me1}}   
\begin{equation}\label{con9}  
{2p-1\choose p-1}\equiv 1 -2p\sum_{k=1}^{p-1}\frac{1}{k}
-2p^2\sum_{k=1}^{p-1}\frac{1}{k^2}\pmod{p^6}. 
    \end{equation}

{\it Remark} 2. Clearly, both congruences 
(\ref{con8}) and (\ref{con9}) can be considered as generalizations
of (\ref{con4}) modulo $p^6$.

Quite recently, in 2011 R. Me\v{s}trovi\'c \cite[Theorem 1.1]{me1}
extended the congruence (\ref{con9}); he proved that 
{\it for any prime $p\ge 11$}
   \begin{equation}\label{con10}
{2p-1\choose p-1}\equiv 1 -2p\sum_{k=1}^{p-1}\frac{1}{k}
+4p^2\sum_{1\le i<j\le p-1}\frac{1}{ij}\pmod{p^7}, 
     \end{equation} 
 which by using the shuffle relation also {\it can be written
in terms of two power sums as}
  \begin{equation}\label{con11}
{2p-1\choose p-1}\equiv 1 -2p\sum_{k=1}^{p-1}\frac{1}{k}
+2p^2\left(\left(\sum_{k=1}^{p-1}\frac{1}{k}\right)^2-
\sum_{k=1}^{p-1}\frac{1}{k^2}\right)\pmod{p^7}. 
     \end{equation} 

{\it Remark} 3. Note that the congruences (\ref{con10}) and (\ref{con11})
reduces to the identity, while for  $p=7$ (\ref{con10}) and (\ref{con11})
are satisfied modulo $7^6$.\hfill $\Box$

Quite recently, R. Tauraso \cite{t2} informed the author that using 
a very similar method  to the method applying in \cite{me1}
to prove the above congruence (\ref{con10}), this congruence  can be  
improved to the following result: 
{\it for any prime $p\ge 7$}
   \begin{equation}\label{con10'}\begin{split}
{2p-1\choose p-1}\equiv & 1 +2p\sum_{k=1}^{p-1}\frac{1}{k}
+\frac{2}{3}p^3\sum_{k=1}^{p-1}\frac{1}{k^3}+
2p^2\left(\sum_{k=1}^{p-1}\frac{1}{k}\right)^2\\
&+\frac{2}{5}p^5\sum_{k=1}^{p-1}\frac{1}{k^5}
+\frac{4}{3}p^4\left(\sum_{k=1}^{p-1}\frac{1}{k}\right)
\left(\sum_{k=1}^{p-1}\frac{1}{k^3}\right)
\pmod{p^9}. 
     \end{split}\end{equation} 
Here we noticed that, using the method applied in \cite[Lemmas 2.2-2.4]{me1},
the term $\sum_{k=1}^{p-1}1/k^5$ 
on the right of (\ref{con10'}) can be eliminated  to obtain that 
 {\it for any prime $p\ge 7$}
   \begin{equation}\label{con11'}\begin{split}
{2p-1\choose p-1}\equiv &  1 +p\sum_{k=1}^{p-1}\frac{1}{k}
-\frac{p^2}{2}\left(5\left(\sum_{k=1}^{p-1}\frac{1}{k}\right)^2+
\sum_{k=1}^{p-1}\frac{1}{k^2}\right)\\
&-\frac{p^3}{30}\left(15\left(\sum_{k=1}^{p-1}\frac{1}{k}\right)
\left(\sum_{k=1}^{p-1}\frac{1}{k^2}\right)-2
\sum_{k=1}^{p-1}\frac{1}{k^3}\right)\\
&+\frac{p^4}{40}\left(35\left(\sum_{k=1}^{p-1}\frac{1}{k^2}\right)^2-
26\sum_{k=1}^{p-1}\frac{1}{k^4}\right)
\pmod{p^9}. 
     \end{split}\end{equation} 

{\it Remark} 4. A computation via {\tt Mathematica} 
verifies that both congruences (\ref{con10'}) and (\ref{con11'}) hold.

\section{Wolstenholme's type congruences in terms of Bernoulli numbers}

The {\it Bernoulli numbers} $B_k$ ($k\in\mathbf N$) 
are defined by the generating function
   $$
\sum_{k=0}^{\infty}B_k\frac{x^k}{k!}=\frac{x}{e^x-1}\,.
  $$
It is easy to find the values $B_0=1$, $B_1=-\frac{1}{2}$, 
$B_2=\frac{1}{6}$, $B_4=-\frac{1}{30}$, and $B_n=0$ for odd $n\ge 3$. 
Furthermore, $(-1)^{n-1}B_{2n}>0$ for all $n\ge 1$. 
These and many other properties can be found, for instance, in \cite{ir}
or \cite{gkp}. 

{\it The Glaisher's congruence} (\ref{con5}) {\it involving Bernoulli number
$B_{p-3}$ may be written as 
 \begin{equation}\label{con12}
{2p -1\choose p-1}\equiv 1-\frac{2}{3}p^3B_{p-3}\pmod{p^4}
  \end{equation}
for all primes $p\ge 7$.}

More generaly,  J.W.L. Glaisher  (\cite[p. 21]{gl1}, \cite[p. 323]{gl2})
proved that {\it for any positive integer $n\ge 1$ and 
any prime $p\ge 5$} 
  \begin{equation}\label{con13}
{np -1\choose p-1}\equiv 1-\frac{1}{3}n(n-1)p^3B_{p-3}\pmod{p^4}.
  \end{equation}
{\it Also, the  congruence $(\ref{con6})$ 
$($cf. the congruence $(\ref{con15})$
below) in terms of Bernoulli numbers may be written as 
 \begin{equation}\label{con14}
{2p -1\choose p-1}\equiv 1-p^3B_{p^3-p^2-2}\pmod{p^5}.
  \end{equation}
for each prime $p\ge 7$.}

In 2008 C. Helou and G. Terjanian \cite{ht} established 
many  Wolstenholme's type congruences modulo $p^k$ with a prime $p$
and $k\in\mathbf N$ such that $k\le 6$. 
As an application, by \cite[Corollary 2(2), p. 493 (also see Corollary 
6(2), p. 495)]{ht}), {\it for any prime $p\ge 5$ we have} 
   \begin{equation}\label{con15}  
{2p-1\choose p-1} \equiv 1-p^3B_{p^3-p^2-2}
+\frac{1}{3}p^5B_{p-3}-\frac{6}{5}p^5B_{p-5}\pmod{p^6}.
   \end{equation}
Applying a technique of Helou and Terjanian \cite{ht} 
based on Kummer type congruences, in 2011 R. Me\v{s}trovi\'c 
 \cite[Corollary 1.3]{me1} proved that
{\it the congruence $(\ref{con10})$  
may be expressed in terms of  Bernoulli numbers
as 
  \begin{equation}\label{con16}\begin{split}    
{2p-1\choose p-1} \equiv& 1-p^3B_{p^4-p^3-2}
+p^5\left(\frac{1}{2}B_{p^2-p-4}-2B_{p^4-p^3-4}\right)\\
&+p^6\left(\frac{2}{9}B_{p-3}^2-\frac{1}{3}B_{p-3}-\frac{1}{10}B_{p-5}\right)
\pmod{p^7}
  \end{split}\end{equation}    
for all primes $p\ge 11$.}

{\it Remark} 5. Note that reducing the moduli and using 
the {\it Kummer congruences} presented in \cite{ht}, 
from (\ref{con16}) may be easily deduced the congruence  (\ref{con15}). 
\hfill $\Box$

\section{Wolstenholme's type harmonic series congruences}

Here, as usually in the sequel, we consider the  congruence relation 
modulo a prime power $p^e$ extended to the ring of rational numbers
with denominators not divisible by $p$. 
For such fractions we put $m/n\equiv r/s \,(\bmod{\,\,p^e})$ 
if and only if $ms\equiv nr\,(\bmod{\,\,p^e})$, and the residue
class of $m/n$ is the residue class of $mn'$ where 
$n'$ is the inverse of $n$ modulo $p^e$.  

As noticed in Section 2, in 1862 J. Wolstenholme \cite{w}
proved that {\it for any prime $p\ge 5$}
 \begin{equation}\label{con17}  
1+\frac{1}{2}+ 
\frac{1}{3}+\cdots +\frac{1}{p-1}\equiv 0\pmod{p^2}.
 \end{equation}
This is in fact an equivalent reformulation of 
Wolstenholme's theorem given by the  congruence (\ref{con4}). 
Wolstenholme \cite{w} also proved that {\it for any prime $p\ge 5$}
 \begin{equation}\label{con18}  
1+\frac{1}{2^2}+ \frac{1}{3^2}+\cdots +\frac{1}{(p-1)^2}\equiv 0\pmod{p}.
 \end{equation}

E. Alkan \cite[Theorem 2, p. 1001]{a} in 1994 proved that {\it for each 
prime $p\ge 5$ the numerator of the fraction
  $$
\frac{1}{1(p-1)}+\frac{1}{2(p-2)}+\cdots +
\frac{1}{\left(\frac{p-1}{2}\right)\left(\frac{p+1}{2}\right)}
  $$
is divisible by $p$} and also noticed that the congruence 
(\ref{con17})  can be deduced from it.

{\it  Remark} 6. In 1999 W. Kimball and W. Webb \cite{kw}
 established the analogue of the congruence (\ref{con17}) in terms of 
Lucas sequences which in particular case reduces to (\ref{con17}). 
Their result is generalized in 2008  by H. Pan \cite[Theorem 1.1]{pan2}.
 \hfill $\Box$

Generalizations of (\ref{con17}) and (\ref{con18}) were established 
by M. Bayat \cite[Theorem 3]{b} in 1997 
(also cf. \cite{ge1} and \cite[Lemma 2.2 and Remark 2.3]{z1}) as follows. 
{\it If $m$ is a positive 
integer and  $p$  a prime such that $p\ge m+3$, then}
  \begin{equation}\label{con19}\begin{split}    
\sum_{k=1}^{p-1}\frac{1}{k^m}\equiv\left\{
    \begin{array}{ll}
0 & \pmod{p}\quad if\,\, m\,\, is\,\, even\\
0 & \pmod{p^2}\quad if\,\, m\,\, is\,\, odd.
  \end{array}\right.
         \end{split}\end{equation}

{\it Remark} 7. For $j=1,2,3$ the numerators of {\it harmonic numbers}
    $$
H_j(n):=\sum_{k=1}^{n}\frac{1}{k^j},\,\, n=1,2,3,\ldots
    $$
written in reduced form are Sloane's sequences 
\cite[sequences A001008, A007406 and A007408]{slo}, respectively.\hfill $\Box$

{\it Remark} 8. 
For a given prime $p$, in \cite{el} and \cite{bo} the authors considered
and investigated  the set $J(p)$ of $n$ for which $p$ 
divides the numerator  of the harmonic sum $H_n:=\sum_{k=1}^{p-1}1/k$.
It is conjectured in \cite[Conjecture A on page 250]{el} 
that the set $J(p)$ is finite for all primes $p$.
This conjecture is recently generalized by J. Zhao \cite{z3}.
\hfill $\Box$

In 1900 J.W.L. Glaisher (\cite[pp. 333-337]{gl2}; also see  
\cite[(v) and (vi) on page 271]{gl3}) proved the following 
generalizations of the congruences (\ref{con18}) and (\ref{con19}) 
(also see \cite[Theorem 5.1(a) and Corollary  5.1]{s1}).
{\it if $m$ is a positive integer and  $p$  a prime such that $p\ge m+3$, then}
  \begin{equation}\label{con20}\begin{split}    
\sum_{k=1}^{p-1}\frac{1}{k^m}\equiv\left\{
    \begin{array}{ll}
\frac{m}{m+1}pB_{p-1-m} & \pmod{p^2}\quad if\,\, m\,\, is\,\, even\\
-\frac{m(m+1)}{2(m+2)}p^2B_{p-2-m} & \pmod{p^3}\quad if\,\, m\,\, is\,\, 
odd.  \end{array}\right.
         \end{split}\end{equation}

{\it Remark} 9. In 2000 Z.H. Sun \cite[Remark 5.1]{s1}
established generalizations modulo $p^4$  of both parts of the 
congruence (\ref{con20}). \hfill $\Box$

In particular, taking $m=1,2,3$ into  the congruence (\ref{con20}) 
we obtain  that {\it for each prime $p\ge 5$} 
  \begin{equation}\label{con21}
\sum_{k=1}^{p-1}\frac{1}{k}\equiv -\frac{1}{3}p^2B_{p-3}\pmod{p^3},
   \end{equation} 
\begin{equation}\label{con22}
\sum_{k=1}^{p-1}\frac{1}{k^2}\equiv \frac{2}{3}pB_{p-3}\pmod{p^2}
   \end{equation} 
and that {\it for each prime $p\ge 7$} 
  \begin{equation}\label{con23}
\sum_{k=1}^{p-1}\frac{1}{k^3}\equiv -\frac{6}{5}p^2B_{p-5}\pmod{p^3}.
   \end{equation} 

{\it Remark} 10. The congruence (\ref{con20}) was also proved in 1938 
by E. Lehmer \cite{leh}.  This congruence 
was generalized in 2007 by X. Zhou and T. Cai \cite[Lemma 3]{zc}
to  multiple harmonic sums; also
see \cite[Theorem 2.14]{z1}.\hfill $\Box$

Another generalization of the congruence (\ref{con17})
is due in 1954 by L. Carlitz \cite{ca1}: 
{\it if $m$ is an arbitrary integer, 
then for each prime $p\ge 5$}
\begin{equation}\label{con24}  
\frac{1}{mp+1}+ \frac{1}{mp+2}+ 
\frac{1}{mp+3}+\cdots +\frac{1}{mp+(p-1)}\equiv 0\pmod{p^2}.
 \end{equation}

{\it Remark} 11. The congruence (\ref{con24}) was also proved in 1989 by 
S. Zhang \cite{zhan}.\hfill $\Box$

Using $p$-adic $L$-functions,
that is the Washington's $p$-adic expansion of the sum
$\sum_{k=1\, (k,p)=1}^{np}1/k^r$ \cite{wa}, in 2000  
S. Hong \cite[Theorem 1.1]{hon2} proved the following generalization of a 
Glaisher's  congruence
(\ref{con20}).
{\it Let $p$ be an odd prime and let $n\ge 1$ and $r\ge 1$
be integers. Then}
 \begin{equation}\label{con25}\begin{split}    
\sum_{k=1}^{p-1}\frac{1}{(np+k)^r}\equiv\left\{
    \begin{array}{ll}
-\frac{(2n+1)r(r+1)}{2(r+2)}p^2B_{p-r-2} & \pmod{p^3}
\, if\, r\, is\, odd\\
& and\, p\ge r+4\\
\frac{r}{r+1}pB_{p-r-1} & \pmod{p^2}
\, if\,r\, is\, even\\
& and\, p\ge r+3\\
-(2n+1)p & \pmod{p^2}
\, if\,\, r=p-2.\\
  \end{array}\right.
       \end{split}\end{equation}

  In 2002 Slavutskii \cite{sl3} 
showed how a more general sums (i.e., the sum (\ref{con25}) with a 
power $p^t,t\in\Bbb N$, instead of $p$)
may be studied by elementary methods without the help of $p$-adic 
$L$-functions. Namely, by \cite[Theorem 1.2]{sl3} 
{\it if $p$ is an odd prime, $n\ge 1$,  $r\ge 1$ and $l\ge 1$ 
are integers, then}
 \begin{equation}\label{con25'}\begin{split}    
\sum_{k=1\atop (k,p)=1}^{p^l-1}\frac{1}{(np^l+k)^r}\equiv\left\{
    \begin{array}{ll}
-\frac{(2n+1)r(r+1)}{2(p^{l-1}+r+1)}p^{2l}B_{p^l-p^{l-1}-r-1} &\pmod{p^{3l}}
\, if\,\, r \\
& is\, odd\, and\, p\ge r+4\\
\frac{r}{p^{2l-2}+r}p^lB_{p^{2l-1}-p^{2l-2}-r} & \pmod{p^{3l-1}}
\, if\,r\\
& is\, even\, and\, p\ge r+3\\
-(2n+1)p^{2l-1} & \pmod{p^{2l}}
\, if\\
& r=p-2.\\
  \end{array}\right.
       \end{split}\end{equation}

{\it Remark} 12. It is obvious that taking 
$l=1$ into (\ref{con25'}), we immediately obtain (\ref{con25}). 
In \cite{me5} R. Me\v{s}trovi\'c proved the congruence (\ref{con25}) by using 
very simple and elementary number theory method.\hfill $\Box$

\section{Wolstenholme's type supercongruences}

A. Granwille \cite{gr} established broader generalizations of 
Wolstenholme's theorem. As an application, it is obtained in \cite{gr}
that {\it for a prime $p\ge 5$  there holds}
    \begin{equation}\label{con26}
{2p-1\choose p-1}\Big/{2p\choose p}^3\equiv
{3\choose 2}\Big/{2\choose 1}^3\pmod{p^5}.
   \end{equation}
   Moreover, by studying  Fleck's quotients,
in 2007 Z.W. Sun and D. Wan \cite[Corollary 1.5]{sw} discovered a new extension of 
Wolstenholme's congruences. In particular, their result 
yields Wolstenholme's theorem and  
{\it for a prime $p\ge 7$ the following new curious congruence}
   \begin{equation}\label{con27}
{4p-1\choose 2p-1} \equiv {4p\choose p}-1 \pmod{p^5}.
   \end{equation}
Further, the congruence (\ref{con34}) given in the next section 
(also cf. the congruences \cite[7.1.10 and 7.1.11]{no})
immediately implies {\it that for a prime $p\ge 5$
  \begin{equation}\label{con28}
{2p^2\choose p^2} \equiv {2p\choose p} \pmod{p^6}
   \end{equation}
and 
  \begin{equation}\label{con29}
{2p^3\choose p^3} \equiv {2p^2\choose p^2} \pmod{p^9}.
   \end{equation}}
{\it If $p$ is a prime, $k$, $n$ and $m$  are  positive integers such that 
$n\ge m$, ${n\choose m}$ is not divisible by $p$  and 
$m\equiv n\,(\bmod{\,p^k})$, then \cite[7.1.16]{no}
  \begin{equation}\label{con30}
{n\choose m} \equiv {\left[n/p\right]\choose \left[m/p\right]} \pmod{p^k},
   \end{equation}
where $[a]$ denotes the integer part of a real number $a$.}

A harmonic type Wolstenholme's supercongruence
is established in 2011 by Y. Su, J. Yang and S. Li in 
\cite[Theorem, p. 500]{syl}
as follows: {\it if $m$ is any integer, $r$ is a non-negative 
integer and $p\ge 5$ is a prime such that $p^r\mid 2m+1$, then}
  \begin{equation}\label{con31}
\sum_{k=1}^{p-1}\frac{1}{mp+k}\equiv 0\pmod{p^{r+2}}.
  \end{equation}

{\it Remark} 13. Note that the congruence (\ref{con31})
with $r=0$ reduces to the congruence (\ref{con24}). \hfill $\Box$

\section{Ljunggren's and  Jacobsthal's binomial congruences}

By Glaisher's congruence $(\ref{con13})$ (\cite[p. 21]{gl1}, 
\cite[p. 323]{gl2}),   
{\it for any positive integer $n$ and a prime $p\ge 5$} 
    \begin{equation*}
{np-1\choose p-1}\equiv 1\pmod{p^3},
      \end{equation*}
which by the identity ${np\choose p}=n{np-1\choose p-1}$
yields \cite[the congruence 7.1.5]{no}  
  \begin{equation}\label{con32}
{np\choose p}\equiv n\pmod{p^3}.
      \end{equation}
In 1952 Ljunggren  generalized the above congruence as follows 
(\cite{bs}; also see \cite[Theorem 4]{ba}, \cite{gr} and 
\cite[Problem 1.6 (d)]{st}):
{\it if $p\ge 5$ is a prime, $n$ and $m$ are positive integers
with $m\le n$, then}
   \begin{equation}\label{con33}
{np\choose mp} \equiv {n\choose m} \pmod{p^3}.
   \end{equation}

{\it Remark} 14. Note that the  congruence (\ref{con33}) with $m=1$ and 
$n=2$ reduces to the Wolstenholme's congruence 
(\ref{con4}). Furthermore, the combinatorial proof 
of (\ref{con33}) regarding  modulo $p^2$ can be found 
in \cite[Exercise 14(c) on page 118]{st}.\hfill $\Box$

Further, the congruence (\ref{con33})
is refined in 1952 by E. Jacobsthal  (\cite{bs}; also see \cite{gr}
and \cite[Section 11.6, p. 380]{co})
as follows: {\it if $p\ge 5$ is a prime, $n$ and $m$ are positive integers
with $m\le n$, then
    \begin{equation}\label{con34}
{np\choose mp} \equiv {n\choose m} \pmod{p^t},
   \end{equation}
where $t$ is the power of $p$ dividing 
$p^3nm(n-m)$ $($this exponent $t$ can only be 
increased if $p$ divides  $B_{p-3}$, the $(p-3)$rd  Bernoulli number$)$}.

{\it Remark} 15. In the literature, the congruence 
(\ref{con34}) is often called {\it Jacobsthal-Kazandzidis congruence}
(see e.g., \cite[Section 11.6, p. 380]{co}).\hfill $\Box$

In particular, the congruence (\ref{con34}) implies that 
{\it for all nonnegative integers $n$, $m$, $a$, $b$ and $c$ with 
$c\le b\le a$, and any prime $p\ge 5$} 
 \begin{equation}\label{con35}
{np^a\choose mp^b} \equiv {np^{a-c}\choose mp^{b-c}}  \pmod{p^{3+a+2b-3c}}.
   \end{equation} 
Moreover, taking  $a=b$ and $c=1$ into (\ref{con35}) 
(cf. \cite[Section 2, Lemma A]{ge2}, or for a direct proof see
\cite[Lemma 19 of Appendix]{gru}), we find that {\it for any prime $p\ge 5$}
   \begin{equation}\label{con36}
{np^a\choose mp^a} \equiv {np^{a-1}\choose mp^{a-1}}  \pmod{p^{3a}}.
   \end{equation} 

Using elementary method, in 1988 N. Robbins \cite[Theorem 2.1]{rob} 
proved the following result. {\it Let $p\ge 3$ be a prime and let  
$n,m,a,b$ be nonnegative integers with $0\le b\le a$, $0<m<np^{a-b}$
and $nm\not\equiv 0\,(\bmod{\,p})$. Then}
  \begin{equation}\label{con36'}
{np^a\choose mp^b} \equiv {np^{a-b}\choose m}  \pmod{p^a}.
   \end{equation}

{\it Remark} 16. Because the original source \cite{bs}
is not easily accessible, the congruence (\ref{con34})
was rediscovered by various authors, 
including  G.S. Kazandzidis (\cite{kaz} and \cite{kaz2})  in 1968
(its proof is based on $p$-adic method)  and  Yu.A. Trakhtman \cite{tr} in 1974.
Furthermore, in 1995 A. Robert and M. Zuber \cite{rz} 
(see also \cite[Chapter 7, Section 1.6]{ro})
proposed a simple proof of the congruence (\ref{con34}) based 
on well-known properties of the $p$-adic Morita gamma function 
$\Gamma_p$.\hfill $\Box$ 

In 2008 Helou and Terjanian \cite[(1) of Corollary on page 490]{ht} 
proved that {\it if $p\ge 5$ is a prime, $n$ and $m$ are positive integers
with $m\le n$, then
    \begin{equation}\label{con37}
{np\choose mp} \equiv {n\choose m} \pmod{p^s},
   \end{equation}
where $s$ is the power of $p$ dividing $p^3m(n-m){n\choose m}$.}

{\it Remark} 17. It is pointed out in \cite[Remark 6 on page 490]{ht} 
that for a prime $p\ge 5$ using $p$-adic methods, the modulus $p^s$ in the 
congruence (\ref{con37}) can be improved to $p^f$, where $f$ is the power of 
$p$ dividing $p^3mn(n-m){n\choose m}$. Clearly, this result  would be
an improvement of Jacobsthal-Kazandzidis congruence given by (\ref{con34}). 
Notice also that  Z.W. Sun and D.M. Davis \cite[Lemma 3.2]{sd}
via elementary method proved the conguence (\ref{con34}) modulo 
$p^s$, where $p\ge 3$ is a prime and $s$ is the power of $p$ dividing 
$p^2n^2$.
\hfill $\Box$

The Jacobsthal's congruence (\ref{con34}) is refined in 2007 by J. Zhao
\cite[Theorem 3.2]{z1} as follows. {\it For a prime $p\ge 7$
define $w_p<p^2$ to be the unique nonnegative integer 
such that $w_p\equiv\left(\sum_{k=1}^{p-1}1/k\right)/p^2\,(\bmod{\,p^2})$.
Then for all positive integers $n$ and $m$ with $n\ge m$}
  \begin{equation}\label{con38}
{np\choose mp}\Big/{n\choose m}\equiv
1+w_pnm(n-m)p^3\pmod{p^5}.
   \end{equation}

{\it Remark} 18. Since $\frac 12 {2p\choose p}={2p-1\choose p-1}$,
taking $n=2$ and $m=1$ into (\ref{con38}), it becomes
  $$
{2p-1\choose p-1}\equiv
1+2w_pp^3\pmod{p^5},
  $$
which is actually (\ref{con7}).\hfill $\Box$

Further, in  2008  Helou and Terjanian \cite[Proposition 2 (1)]{ht}
generalized the congruence (\ref{con38}) modulo $p^6$ 
in the form involving the Bernoulli numbers 
as follows. {\it Let  $p\ge 5$ be a prime.
Then for all positive integers $n$ and $m$ with $n\ge m$}
  \begin{equation}\label{con39}\begin{split}
{np\choose mp}\Big/{n\choose m}
\equiv & 1-mn(n-m)\left(\frac{p^3}{2}B_{p^3-p^2-2}-
\frac{p^5}{6}B_{p-3}\right.\\
&\left. +\frac{1}{5}(m^2-mn+n^2)p^5B_{p-5}\right)\pmod{p^6}.
   \end{split}\end{equation}
As an application of (\ref{con39}), it can 
be obtained \cite[the congruence of Remark 4 on page 489]{ht}
that for {\it each prime $p\ge 5$}
  \begin{equation}\label{con40}
{np\choose mp}\Big/{n\choose m}\equiv
1-\frac{1}{3}mn(n-m)p^3B_{p-3}\pmod{p^4}.
   \end{equation}

{\it Remark} 19. Taking $m=1$ into (\ref{con40}), it immediately
reduces to the Glaisher's congruence (\ref{con13})
which for $n=2$ becomes (\ref{con12}).
\hfill $\Box$

{\it Remark} 20. In \cite[Sections 3 and 4]{ht} 
C. Helou and G. Terjanian  established 
numerous Wolstenholme's type congruences 
of the form ${np\choose mp}\equiv {n\choose m}P(n,m,p)
\,\,(\bmod\,\, p^k)$ and ${mnp\choose np}\equiv  {mn\choose n}
P(n,m,p)\,\,(\bmod\,\, p^k)$,
where $p$ is  a prime, $k\in\{4,5,6\}$,  $m,n,\in \Bbb N$ with 
$m\le n$, and $P$ is a polynomial of $m,n$ and $p$ involving 
Bernoulli numbers as its coefficients.\hfill $\Box$

{\it Remark} 21. In 2011 R. Me\v{s}trovi\'c 
\cite{me3}  discussed  the following type congruences: 
    $$
{np^k\choose mp^k} \equiv {n\choose m} (\bmod{\,p^r}),
     $$
where $p$ is a prime, $n,m,k$ and $r$ are various positive integers.
\hfill $\Box$

  \section{Wolstenholme primes}
A prime $p$ is said to be a {\it Wolstenholme prime} 
(see \cite[p. 385]{m} or \cite{me2}; this is the Sloane's sequence 
A088164 from \cite{slo}) 
if it satisfies the congruence 
 $$
{2p-1\choose p-1} \equiv 1 \pmod{p^4}.
 $$
Clearly, this is equivalent to the fact that the 
{\it Wolstenholme quotient} $W_p$ defined as 
  $$
W_p=\frac{{2p-1\choose p-1}-1}{p^3}, \,\, p\ge 5
  $$
is divisible by $p$ ($W_p$ is the Sloane's sequence A034602 from \cite{slo}; 
also cf. a related sequence A177783).

Notice that by a special case of 
Glaisher's congruence (\cite[p. 21]{gl1}, \cite[p. 323]{gl2};   
also cf. \cite[Corollary, p. 386]{m}) given by (\ref{con12})
it follows that 
 \begin{equation}\label{con41}
W_p\equiv -\frac{2}{3}B_{p-3}\pmod{p}.
  \end{equation}
From the congruence (\ref{con41}) we see that 
$p$ {\it is a Wolstenholme prime if and only if
$p$ divides the numerator of the Bernoulli number $B_{p-3}$}.
This is by (\ref{con5}) also equivalent with the fact (e.g., see \cite{gar} 
or \cite[Theorem 2.8]{z1}) that the numerator of the fraction 
$\sum_{k=1}^{p-1}1/k$ written in reduced form is  divisible by $p^3$, and 
also by (\ref{con6}) with the fact that the numerator of the fraction 
$\sum_{k=1}^{p-1}1/k^2$ written in reduced form is divisible by $p^2$.
 
In other words, a   Wolstenholme prime is a prime 
$p$ such that $(p,p-3)$ is an {\it irregular pair} 
(see \cite{jo} and \cite{bcem}).
The Wolstenholme primes therefore form a subset of the 
irregular primes (see e.g., \cite[p. 387]{m}). 

In 1995 Mcintosh \cite[Proof of Theorem 2]{m} observed that
by a particular case of a result of Stafford and Vandiver \cite{sv} in
1930 and Fermat's Little Theorem, {\it for any prime $p\ge 11$}
  \begin{equation}\label{con42}
B_{p-3}\equiv \frac{1}{21}\sum_{k=\left[p/6\right]+1}^{\left[p/4\right]}
\frac{1}{k^3}\pmod{p}.
   \end{equation}

{\it Remark} 22. The congruence (\ref{con42}) is very useful 
in the computer search of Wolstenholme primes (see \cite{mr}). 
\hfill $\Box$

Only  two Wolstenholme primes are known today:
16843 and 2124679. The first was found
(though not explicitly reported) by Selfridge and Pollack in 1964
(Notices Amer. Math. Soc. 11 (1964), 97),  and later 
confirmed by W. Johnson \cite{jo} and S.S. Wagstaff 
(Notices Amer. Math. Soc. 23 (1976), A-53).
The second was discovered by Buhler, Crandall,
Ernvall and Mets\"{a}nkyl\"{a} in 1993 \cite{bcem}. 
In 1995, McIntosh \cite{m} determined by calculation that there is no other 
Wolstenholme prime $p<5\cdot 10^8$.
In 2007 R.J. McIntosh and E.L. Roettger \cite{mr} reported  that these
primes are only two  Wolstenholme primes less than $10^9$.
However, using the argument based on the prime number theorem, in 1995
McIntosh  \cite[p. 387]{m} conjectured that there are infinitely many 
Wolstenholme primes. It seems that the proof of this assertion 
would be very difficult.

In 2007 J. Zhao \cite{z2} defined a Wolstenholme prime via harmonic numbers;
namely,  a prime $p$ is a {\it Wolstenholme prime if} 
   $$
\sum_{k=1}^{p-1}\frac{1}{k}\equiv 0 \pmod{p^3}.
   $$
In 2011 R. Me\v{s}trovi\'c \cite[Proposition 1]{me2}
proved that if  {\it $p$ is a Wolstenholme prime, then}
   \begin{equation}\label{con43}\begin{split}  
{2p-1\choose p-1} 
\equiv &1 +p \sum_{k=1}^{p-1}\frac{1}{k}
-\frac{p^2}{2} \sum_{k=1}^{p-1}\frac{1}{k^2}+
\frac{p^3}{3} \sum_{k=1}^{p-1}\frac{1}{k^3}-
\frac{p^4}{4} \sum_{k=1}^{p-1}\frac{1}{k^4}\\
&+\frac{p^5}{5} \sum_{k=1}^{p-1}\frac{1}{k^5}-
\frac{p^6}{6} \sum_{k=1}^{p-1}\frac{1}{k^6}\pmod{p^8}.
   \end{split}\end{equation}
The above congruence can be simplified as follows 
\cite[Proposition 2]{me2}:
  \begin{equation}\label{con44}
{2p-1\choose p-1} \equiv 1 +\frac{3p}{2} \sum_{k=1}^{p-1}\frac{1}{k}
-\frac{p^2}{4}\sum_{k=1}^{p-1}\frac{1}{k^2}
+\frac{7p^3}{12} \sum_{k=1}^{p-1}\frac{1}{k^3}
+\frac{5p^5}{12} \sum_{k=1}^{p-1}\frac{1}{k^5}\pmod{p^8}.
   \end{equation}  
Reducing the modulus in the previous congruence,
we can obtain the following simpler congruences 
for Wolstenholme prime $p$ \cite[Corollary 1]{me2}:
   \begin{equation}\label{con45}\begin{split}    
{2p-1\choose p-1} &\equiv 1 -2p \sum_{k=1}^{p-1}\frac{1}{k}
-2p^2\sum_{k=1}^{p-1}\frac{1}{k^2}\pmod{p^7}\\
&\equiv 1 +2p \sum_{k=1}^{p-1}\frac{1}{k}+
\frac{2p^3}{3}\sum_{k=1}^{p-1}\frac{1}{k^3}\pmod{p^7}.
   \end{split}\end{equation}

In terms of  the Bernoulli numbers the congruence  (\ref{con45}) 
may be written as \cite[Corollary 2]{me2}
 \begin{equation}\label{con46}    
{2p-1\choose p-1} \equiv 1-p^3B_{p^4-p^3-2}
-\frac{3}{2}p^5B_{p^2-p-4}+\frac{3}{10}p^6B_{p-5}\pmod{p^7}.
  \end{equation}
The congruence (\ref{con46}) can be given by the following 
expression involving lower order Bernoulli numbers 
\cite[Corollary 2]{me2}: 
       \begin{equation}\label{con47}\begin{split}  
&{2p-1\choose p-1} \equiv 1-p^3\left(\frac{8}{3}B_{p-3}
-3B_{2p-4}+\frac{8}{5}B_{3p-5}-\frac{1}{3}B_{4p-6}\right)\\
&-p^4\left(\frac{8}{9}B_{p-3}-\frac{3}{2}B_{2p-4}+\frac{24}{25}B_{3p-5}
-\frac{2}{9}B_{4p-6}\right)\\
&-p^5\left(\frac{8}{27}B_{p-3}-\frac{3}{4}B_{2p-4}
+\frac{72}{125}B_{3p-5}-\frac{4}{27}B_{4p-6}+\frac{12}{5}B_{p-5}\right.\\
&\left.-B_{2p-6}\right)-\frac{2}{25}p^6B_{p-5}\pmod{p^7}.
 \end{split}\end{equation}    

{\it Remark} 23. 
Combining the first congruence of (\ref{con45}) 
and a recent result of the author in \cite[Theorem 1.1]{me1}
given by the congruence (\ref{con11}), 
we obtain a new characterization of Wolstenholme primes 
as follows \cite[Remark 1.6]{me1}. \hfill $\Box$

{\it A prime  $p$ is  a Wolstenholme prime if and only if}
  $$  
{2p-1\choose p-1} \equiv 1 -2p \sum_{k=1}^{p-1}\frac{1}{k}
-2p^2\sum_{k=1}^{p-1}\frac{1}{k^2}
 \pmod{p^7}.
   $$
 
{\it Remark} 24 (\cite[Remark 1]{me1}).  
A computation via {\tt Mathematica} 
shows that no prime  $p<10^5$
satisfies the second congruence from (\ref{con45}),
except the Wolstenholme prime $16843$.
Accordingly, an interesting question is as follows:
{\it Is it  true that the second congruence  from $(\ref{con45})$
yields that a prime $p$ is necessarily a Wolstenholme prime}?
{\it We conjecture that this is true.}\hfill $\Box$

\section{Wolstenholme's type theorems for composite moduli}

For positive integers $k$ and $n$ let $(k,n)$ denotes 
the greatest common divisor of $k$ and $n$. 
In 1889 C. Leudesdorf \cite{le} (also see \cite[Ch. VIII]{hw} 
or \cite[Chapter 3, the congruence (15) on page 244]{sc}) was proved 
that {\it for any positive integer $n$ such that $(n,6)=1$
     \begin{equation}\label{con48}
   \sum_{k=1\atop (k,n)=1}^{n-1}\frac{1}{k}\equiv 0\pmod{n^2},
   \end{equation}
where the summation ranges over all $k$ with $(k,n)=1$.}

{\it Remark} 25. Observe that when  $n=p\ge 5$ is a prime, then the congruence 
(\ref{con48}) reduces to the Wolstenholme's congruence 
(\ref{con17}).\hfill $\Box$

{\it Remark} 26. In 1934 S. Chowla \cite{ch} gave a very short and 
elegant proof of the congruence (\ref{con48}). This was in another way also
proved by H. Gupta \cite[the congruence (1)]{gup}.
Furthermore, in 1933 Chowla \cite{ch1} proved  Leudesdorf's theorem when 
$n$ is a power of a prime $p\ge 5$.  \hfill $\Box$

The binomial coefficient's analogue of (\ref{con48}) was established in 1995
by R.J. McIntosh. Namely, in \cite[Section 2]{m}
for positive inegers $n$, the author 
defined the {\it modified binomial coefficients}
 $$
{2n-1\choose n-1}'=\prod_{k=1\atop (k,n)=1}^n\frac{2n-k}{k}
  $$
and observed that for primes $p$, ${2p-1\choose p-1}'={2p-1\choose p-1}$.

{\it Then by \cite[Theorem 1]{m}, for $n\ge 3$
 \begin{equation}\label{con49}
{2n-1\choose n-1}'\equiv 1+n^2\varepsilon_n\pmod{n^3}, 
   \end{equation}
where 
    \begin{equation*}\begin{split}    
\varepsilon_n=\left\{
    \begin{array}{ll}
n/2 & \quad if\,\, n\,\, is\,\, a\,\, power\,\, of \,\,2,\\
(-1)^{r+1}n/3 & \quad if\,\, n\equiv 0\,(\bmod \, 3)\,\, and\,\, 
n\,\, has\,\, exactly\\
&\quad r \,\,distinct\,\, prime\,\, factors,\,\, each\,\, 
\not\equiv 1\,(\bmod \, 6),\\
0 &\quad  otherwise.  \end{array}\right.
         \end{split}\end{equation*}}

{\it Remark} 27. McIntosh \cite[page 384]{m} also noticed that in 
1941 H.W. Brinkman \cite{se} in his partial solution to David Segal's
conjecture observed the following relation 
between the ordinary binomial coefficient and the 
modified binomial coefficient:
    $$
{2n-1\choose n-1}=\prod_{d\mid n}{2d-1\choose d-1}',
  $$
which for example, for $n=p^2$ with a prime $p$ becomes
 $$
{2p^2-1\choose p^2-1}={2p^2-1\choose p-1}'{2p-1\choose p-1}.
  $$

{\it Remark} 28. Some extensions of (\ref{con48}) can be found
in the monograph of G.H. Hardy and E.M. Wright \cite[Ch. VIII]{hw},
\cite{hw2} and \cite{ra}.\hfill $\Box$

Further generalizations of Leudesdorf's congruence 
(\ref{con48}) were obtained by L. Carlitz \cite{ca1} in 1954 and 
by H.J.A. Duparc and W. Peremans \cite[Theorem 2]{dp} in 1955. Their result 
(i.e., \cite[Theorem 2]{dp}) was also established in 
1982 by I. Gessel \cite[Theorems 1 and 2]{ge}
and by I.Sh. Slavutskii  \cite[Corollary 1]{sl1} in 1999. 

Namely, by \cite[Corollary 1(a)]{sl1}, {\it if $s$ is an even positive 
integer and $n$ a positive integer such that $(p,n)=1$ 
for all primes $p$ such that $(p-1)\mid s$, then 
 \begin{equation}\label{con50}
   \sum_{k=1\atop (k,n)=1}^{n-1}\frac{1}{k^s}\equiv 0\pmod{n}.
   \end{equation}}
On the other hand, by \cite[Corollary 1(b)]{sl1}, 
{\it if $s$ is an odd positive 
integer and $n$ a positive integer such that $1)$
$p-1$ don't divide $s+1$ for every prime $p$ with $p\mid n$;
or $2)$ $p\mid s$ for all primes $p$ such that $(p-1)\mid (s+1)$
and $p\mid n$, then}
  \begin{equation}\label{con51}
   \sum_{k=1\atop (k,n)=1}^{n-1}\frac{1}{k^s}\equiv 0\pmod{n^2}.
   \end{equation}
Both congruences (\ref{con50})  and 
(\ref{con51}) are immediate consequences of Theorem 2 in \cite{sl1}
(also see \cite[Theorem 1]{sl2} and \cite[Chapter 3, 
the congruence (15') on page 244]{sc}):
{\it if $s$ and $n$ are positive integers such that 
$(n,6)=1$ and $t=(\varphi(n^2)-1)s$,
where $\varphi(\cdot )$ denotes the Euler's totient function, then}
 \begin{equation}\label{con52}\begin{split}    
\sum_{k=1\atop (k,n)=1}^{n-1}\frac{1}{k^s}\equiv
\left\{
    \begin{array}{ll}
n\prod\limits_{p\mid n}(1-p^{t-1})B_t\qquad\,\,\, (\bmod{\,n^2}) 
& \quad for\,\,even\,\, s\\
\frac{t}{2}n^2\prod\limits_{p\mid n}(1-p^{t-2})B_{t-1} \,\,\, (\bmod{\,n^2})& \quad 
for\,\,odd\,\, s.
  \end{array}\right.
         \end{split}\end{equation}

Notice that  by the  congruence (\ref{con52}) (cf. (6) and (7) in \cite{sl1}), 
it follows that {\it for $n=p^l$ with $l\in\Bbb N$ and a prime $p\ge 5$}
   \begin{equation}\label{con53}\begin{split}    
\sum_{k=1\atop (k,p)=1}^{p^l-1}\frac{1}{k^s}\equiv
\left\{
    \begin{array}{ll}
p^lB_t\,\,\, (\bmod{\,p^{2l}})& \quad for\,\, even\,\, s, t=(\varphi(p^{2l})-1)s\\
\frac{t}{2}p^{2l}B_{t-1}\,\,\, (\bmod{\,p^{2l}})
 & \quad for\,\, odd \,\, s, t=(\varphi(p^{2l})-1)s. \end{array}\right.
    \end{split}\end{equation}
As an application of the congruence (\ref{con53}), we obtain 
the following result obtained in 1955 by H.J.A. Duparc and W. Peremans 
\cite[Theorem 1]{dp} (cf. \cite[Corollary 2]{sl1}). 
{\it Let $n=p^l$ be a power of a prime $p\ge 3$.  Then}
   \begin{equation}\label{con54}\begin{split}    
\sum_{k=1\atop (k,p)=1}^{p^l-1}\frac{1}{k^s}\equiv
\left\{
    \begin{array}{ll}
 0\,(\bmod{p^{2l-1}})&\, for\,odd\,\, s\,\,with\,\, p-1\mid s+1\, and 
\,s\not\equiv\,0\,(\bmod{p})  \\
 0\,(\bmod{p^{2l}})&\, for\,odd\,\, s\,\,with\,\, s+1\not\equiv\,0
\,(\bmod{\,p-1})\,or \, p\mid s  \\
 0\,(\bmod{\,p^{l-1}})&\, for\,even\,\, s\,\,with\,\,  p-1\mid s \\
 0\,(\bmod{\,p^{l}})&\, for\,even\,\, s\,\,with\,\, s\not\equiv\,0\,
(\bmod{\,p-1}).\end{array}\right.
    \end{split}\end{equation}

{\it Remark} 29. Notice that the second part and the fourth part of the 
congruence (\ref{con54}) under the condition $s\le p-3$ 
were also proved in 1997 by M. Bayat \cite[Theorem 4]{b}.
Moreover, the second part and the fourth part of the 
congruence (\ref{con54}) with $n=ap^l$ ($a,l\in\Bbb N$) instead 
of $p^l$, $s=1$ and $s=2$ were also proved  in 1998 by  D. Berend and J.E. 
Harmse \cite[Proposition 2.2]{bh}.  \hfill $\Box$

Wolstenholme type congruence for 
product of distinct primes was established in 2007 by  S. Hong
\cite{hon}: {\it if $m$ and $n$ are integers with $m\ge 0$, $n\ge 1$,
$<n>:=\{1,\ldots,n\}$, $p_1,\ldots, p_n$ are disctinct 
primes and all greater than $3$, then}
   \begin{equation}\label{con55}
 \sum_{k=1\atop \forall i\in <n>, 
(k,p_i)=1}^{p_1\cdots p_n}\frac{1}{mp_1\cdots p_n+k}\equiv 0
\pmod{(p_1\cdots p_n)^2}.
   \end{equation}

\section{On the converse of Wolstenholme's theorem}
If $n\ge 5$ is a prime, then by Wolstenholme's theorem
  \begin{equation*}
{2n-1\choose n-1}\equiv 1\pmod{n^3}.
  \end{equation*}
{\it Is the converse true}? This question, still unanswered today, 
has been asked by J.P. Jones for many years 
(see \cite[Chapter 2, p. 23]{ri}, \cite{gr} and \cite[B31, p. 131]{gu}).

In 2001 V. Trevisan and K.E. Weber \cite[Theorem 1]{tw} proved that 
{\it if $n$ is an even positive integer, then} 
     $$
{2n-1\choose n-1}\not\equiv 1\pmod{n^3}.
      $$
Following \cite[Chapter 2, p. 23]{ri}, the mentioned 
problem leads naturally to the following concepts and 
questions. Let $n\ge 5$  be odd, and let 
  $$
A(n):={2n-1\choose n-1}.
  $$
For each $k\ge 1$ we may consider the set 
  $$
W_k=\{n\,\,{\rm odd}\,\, , n\ge 5|\, A(n)\equiv 1\,(\bmod{\,n^k})\}.
  $$
Obviously, $W_1\supset W_2\supset W_3\supset W_4\supset \ldots .$
From Wolstenholme's theorem every prime number greater than 3 belongs to 
$W_3$. {\it Jones' question is whether $W_3$ is just the set of these prime 
numbers}.    
  
Notice that the set $W_4$ coincides with the set of all Wolstenholme primes
defined in Section 7. The set of composite integers in $W_2$
contains the squares of Wolstenholme primes.  Mcintosh \cite[p. 385]{m}
conjectured that these sets coincide and verified that this is true up to
$10^9$; the only composite number $n\in W_2$ with $n<10^9$, is 
$n=283686649=16843^2$. Furthermore, 
using the argument based on the prime number theorem, 
McIntosh  (\cite[p. 387]{m}) conjectured that the set $W_5$ is empty;
this means that no prime satisfies
the congruence 
$$
{2p-1\choose p-1} \equiv 1\pmod{p^5}.
 $$

Recall also that in 2010 
K.A. Broughan, F. Luca and I.E. Shparlinski \cite{bls} 
investigated the subset $W_1'$ 
consisting of all composite positive integers $n$
belonging to the set $W_1$.  They proved \cite[Theorem 1]{bls} 
that the set $W_1'$ is of asymptotic density zero. More precisely, if 
$W(x)$ is defined to be the number of composite positive integers
$n\le x$ which satisfy ${2n-1\choose n-1} \equiv 1\,\, (\bmod{\,\,n})$,
then $\lim_{x\to\infty}W(x)/x=0$. McIntosh \cite[p. 385]{m}
reported that the only elements of the set $W_1$ less than $10^9$ that are not
primes nor prime powers are $29\times 937$ and $787\times 2543$, 
but none of these satisfy Wolstenholme's congruence.

{\it Remark} 30. The converse of Wolstenholme's theorem 
for particular classes of composite integers $n$ was  discussed 
and proved in 2001 by Trevisan and Weber \cite{tw}.
Further, in 2008 Helou and Terjanian \cite[Section 5, 
Propositions 5-7]{ht} deduced that this converse holds for many infinite 
families of composite 
integers $n$.\hfill $\Box$

\section{Binomial sums related to Wolstenholme's theorem}

In 2006 M. Chamberland and K. Dilcher \cite{chd1} studied a class of binomial 
sums of the form 
   $$
u_{a,b}^{\varepsilon}(n):=\sum_{k=0}^n(-1)^{\varepsilon k}{n\choose k}^a
{2n\choose k}^b,
  $$
for nonnegative integers $a,b,n$ and $\varepsilon \in\{0,1\}$,
and showed that these sums are closely related to Wolstenholme's theorem.
{\it Namely, they proved \cite[Theorem 3.1]{chd1} that for any 
prime $p\ge 5$ holds
 \begin{equation}\label{con56}
u_{a,b}^{\varepsilon}(p):=\sum_{k=0}^p(-1)^{\varepsilon k}{p\choose k}^a
{2p\choose k}^b\equiv 1+(-1)^{\varepsilon}2^b\pmod{p^3},
  \end{equation}
except when $(\varepsilon, a,b)=(0,0,1)$.}

In a subsequent paper in 2009 M. Chamberland and K. Dilcher \cite{chd2}
studied the above sum for $(\varepsilon, a,b)=(1,1,1)$, that is, 
with the simplified notation, the sum 
 $$
u(n):=\sum_{k=0}^n(-1)^{k}{n\choose k}{2n\choose k}.
  $$
{\it Under this notation, the authors proved \cite[Theorem 2.1]{chd2}
that for all primes $p\ge 5$ and integers $m\ge 1$ we have}
  \begin{equation}\label{con57}
u(mp)\equiv u(m)\pmod{p^3}.
  \end{equation}

In 2002 T.X. Cai and A. Granville \cite[Theorem 6]{cg} proved 
the following result. {\it If $p\ge 5$ is a prime and $n$ a positive integer, 
then
 \begin{equation}\label{con58}\begin{split}    
\sum_{k=0}^{p-1}(-1)^k{p-1\choose k}^n\equiv\left\{
    \begin{array}{ll}
{np-2\choose p-1} & \pmod{p^4}\quad if\,\, n\,\, is\,\, odd\\
2^{n(p-1)} & \pmod{p^3}\quad if\,\, n\,\, is\,\, even
  \end{array}\right.
         \end{split}\end{equation}
and} 
 \begin{equation}\label{con59}\begin{split}    
\sum_{k=0}^{p-1}{p-1\choose k}^n\equiv\left\{
    \begin{array}{ll}
{np-2\choose p-1} & \pmod{p^4}\quad if\,\, n\,\, is\,\, even\\
2^{n(p-1)} & \pmod{p^3}\quad if\,\, n\,\, is\,\, odd.
  \end{array}\right.
         \end{split}\end{equation}

In 2009 H. Pan \cite[Theorem 1.1]{pan3}  
generalized the second parts of congruences (\ref{con58})
and (\ref{con59}) as follows. 
{\it Let $p\ge 3$ be a prime and let $n$ be a positive integer. Then}

   \begin{equation}\label{con60'}
\sum_{k=0}^{p-1}(-1)^{(n-1)k}{p-1\choose k}^a\equiv 2^{n(p-1)}+
\frac{n(n-1)(3n-4)}{48}p^3B_{p-3}\pmod{p^4}.
   \end{equation}

Recently, in 2011 R. Me\v{s}trovi\'c \cite[Theorem 3]{me4}
extended Pan's congruence (\ref{con60'}) for $n=-1$
by proving the following congruence for the sum of the reciprocals of 
binomial coefficients.  {\it Let $p\ge 3$ be a  prime. Then
   \begin{equation}\label{con61}
\sum_{k=0}^{p-1}{p-1\choose k}^{-1}\equiv 2^{1-p}-
\frac{7}{24}p^3B_{p-3}\pmod{p^4}.
   \end{equation}
In particular, we have}    
   \begin{equation}\label{con62}
\sum_{k=0}^{p-1}{p-1\choose k}^{-1}\equiv 2^{1-p}\pmod{p^3}.
    \end{equation}

Many interesting congruences modulo $p^k$ with $k\ge 3$
are in relation to the {\it Ap\'{e}ry numbers} 
$A_n$ defined in 1979 by R. Ap\'ery \cite{ap} as
  $$
A_n:=\sum_{k=0}^n{n\choose k}^2{n+k\choose k}^2=
\sum_{k=0}^n{n+k\choose 2k}^2{2k\choose k}^2,\,n=0,1,2,\ldots.
  $$
For example, in 1982 I. Gessel \cite{ge} 
 proved that {\it for any prime $p\ge 5$}
  \begin{equation}\label{con63}
A_{pn}\equiv A_n\pmod{p^3},\,n=0,1,2,\ldots .
  \end{equation} 

{\it Remark} 31. Z.W. Sun in \cite{su1}, \cite[pp. 48-49]{su2} 
and \cite{su3} also made many interesting 
conjectures on congruences involving the Ap\'ery numbers $A_n$.\hfill $\Box$

Finally, we present here an interesting congruence 
proposed as a problem on W.L. Putnam Mathematical Competition
\cite{pu2}. {\it If $p\ge 5$ is a prime and $k=[2p/3]$, then 
by \cite[Problem A5 (1996)]{pu2}}
  \begin{equation}\label{con64}
{p\choose 1}+{p\choose 2}+\cdots +{p\choose k}\equiv 0\pmod{p^2}.
  \end{equation}

\section{$q$-analogues of Wolstenholme's type congruences}

Recal that the {\it generalized harmonic numbers} $H_n^{(m)}$,
$n,m=0,1,2,\ldots$ are defined by
    $$
H_n^{(m)}=\sum_{k=1}^n\frac{1}{k^m}
  $$
(we assume that $H_0^{(m)}=0$ for all $m$). Notice that
  $$
H_n^{(1)}:=H_n=\sum_{k=1}^n\frac{1}{k} 
  $$
is the {\it harmonic number}. A $q$-{\it analog} of $H_n$ is given by 
the $q$-{\it harmonic numbers}
  $$
H_n(q):=\sum_{k=1}^n\frac{1}{[k]_q},\quad n\ge 0,\,\, |q|<1,
  $$
where 
 $$
[k]_q:=\frac{1-q^k}{1-q}=1+q+\cdots +q^{k-1}.
 $$
A different $q$-{\it analog} of $H_n$ is 
 $$
\widetilde{H}_n(q):=\sum_{k=1}^n\frac{q^k}{[k]_q},\quad n\ge 0,\,\, |q|<1,
  $$
In 1999 G.E. Andrews \cite[Theorem 4]{an} proved a $q$-analogue of the weaker 
version (modulo $p$) of the congruence (\ref{con17}); {\it namely, 
for pimes $p\ge 3$}
   \begin{equation}\label{con65}
H_{p-1}(q)\equiv \frac{p-1}{2}(1-q)\pmod{[p]_q}.
   \end{equation}
Andrews  also proved that  {\it for pimes $p\ge 3$}
   \begin{equation}\label{con66}
\widetilde{H}_{p-1}(q)\equiv -\frac{p-1}{2}(1-q)\pmod{[p]_q}.
   \end{equation}
In 2007 L.L. Shi and H. Pan \cite[Theorem 1]{sp}
(also see \cite[the congruence (1.3)]{pan}) extended (\ref{con65}) to
  \begin{equation}\label{con67}
H_{p-1}(q)\equiv \frac{p-1}{2}(1-q)+\frac{p^2-1}{24}(1-q)^2[p]_q\pmod{[p]_q^2}
   \end{equation}
{\it for each prime $p\ge 5$}. 

In 2007 L.L. Shi and H. Pan \cite[Lemma 2 (5) and (4)]{sp} also showed that
{\it for each prime $p\ge 5$
  \begin{equation}\label{con68}
\sum_{k=1}^{p-1}\frac{1}{[k]_q^2}\equiv 
-\frac{(p-1)(p-5)}{12}(1-q)^2\pmod{[p]_q}.
   \end{equation}
and 
\begin{equation}\label{con69}
\sum_{k=1}^{p-1}\frac{q^k}{[k]_q^2}\equiv 
-\frac{p^2-1}{12}(1-q)^2\pmod{[p]_q}.
   \end{equation}}

Recently, in 2011 A. Straub \cite{str} proved a $q$-analogue
of a classical binomial congruence (\ref{con33}) due to Ljunggren.
If under the above notation we set
   $$
[n]_q!:=[n]_q[n-1]_q\cdots [1]_q
  $$
and 
  $$
{n\choose k}_q:=\frac{[n]_q!}{[k]_q![n-k]_q!}
  $$
(this is a polynomial in $q$ with integer coefficients)
then Straub \cite[Theorem 1]{str} proved that 
{\it if $p\ge 5$ is a prime, $n$ and $m$ are  nonnegative integers
with $m\le n$, then}
   \begin{equation}\label{con70}
{np\choose mp}_q \equiv {n\choose m}_{q^{p^2}}-
{n\choose m+1}{m+1\choose 2}\frac{p^2-1}{12}(q^p-1)^2 \pmod{[p]_q^3}.
   \end{equation} 
Notice that {\it the congruence $(\ref{con70})$ reduced modulo
$[p]_q^2$  becomes} 
 \begin{equation}\label{con71}
{np\choose mp}_q \equiv {n\choose m}_{q^{p^2}}\pmod{[p]_q^2},
   \end{equation} 
which was proved in 1995 by W.E. Clark 
\cite[the congruence (2) on page 197]{cla}. 

Notice also that in 1999 G.E. Andrews \cite{an} proved
a similar result; e.g.:
  \begin{equation}\label{con72}
{np\choose mp}_q \equiv q^{(n-m)m{p\choose 2}}{n\choose m}_{q^p}\pmod{[p]_q^2}.
   \end{equation} 

Furthermore, taking $m=2$ and $n=1$ 
into the congruence (\ref{con70}), we obtain 
{\it a $q$-analogue of Wolstenholme's theorem as}:
  \begin{equation}\label{con73}
{2p\choose p}_q \equiv [2]_{q^{p^2}}-
\frac{p^2-1}{12}(q^p-1)^2 \pmod{[p]_q^3}.
   \end{equation} 

{\it Remark} 32.  The congruences in this section 
are to be understood 
as congruences in the polynomial ring $\Bbb Z[q]$. Note that
it is clear that $[p]_q:=1+q+\cdots +q^{p-1}$, 
as the $p$th {\it cyclotomic polynomial} is irreducible; hence the 
denominator of $H_{p-1}(q)$, seen as a 
rational function of $q$, is relatively prime to $[p]_q$.\hfill $\Box$

{\it Remark} 33. In 2008 K. Dilcher \cite[Theorems 1 and 2]{dil} 
generalized the congruences (\ref{con65}) and (\ref{con66})
deriving congruences $(\bmod{\,[p]_q})$ for the {\it generalized} 
(or {\it higher-order}) $q$-{\it harmonic numbers}. His results are in fact
$q$-analogues of the congruences 
$H_{p-1}^{(k)}\equiv 0 \,(\bmod{\, p})$ which follow from 
(\ref{con19}).\hfill $\Box$

{\it Note.} I found References \cite{che}, 
\cite{sh1}, \cite{sh2}, \cite{ta}, \cite{ya}, \cite{ylz} and \cite{zy} on
the Internet, but they are not cited in this article
because of they are not accessible to the author.

\vspace{4mm}

\centerline{A{\scriptsize{PPENDIX}}}

\vspace{3mm}

\centerline{\bf{A) External Links on Wolstenholme's theorem and 
Wolstenholme primes}}

{\small Eric Weisstein World of Mathematics, Wolstenholme prime, 

\noindent{\tt http://mathworld.wolfram.com/Wolstenholme prime.html}, 
from MathWorld. 

Wikipedia {\tt http://en.wikipedia.org/wiki/Wolstenholme\_prime}

{\tt http://planetmath.org/encyclopedia/WolstenholmesTheorem.html}

{\tt http://imperator.usc.edu/~bruck/research/stirling/}

{\tt http://mathforum.org/kb/thread.jspa?messageID=29014\&tstart=0}

{\tt http://www.music.us/education/W/Wolstenholme's-theorem.htm}

{\tt http://www.music.us/education/W/Wolstenholme's-prime.htm}

{\tt http://www.music.us/education/J/Joseph-Wolstenholme.htm}

{\tt http://chandrumath.wordpress.com/2010/10/03/Wolstenholmes-theorem}

{\tt http://mathoverflow.net/question/26137/binomial}

\noindent{\tt -supercongrueces-is-there-any-reason-for-them}

{\tt http://uniblogger.com/en/Wolstenholme's$\_$theorem}}

\vspace{3mm}

\centerline{\bf{B) Sloane's sequences related to Wolstenholme's
theorem and Wolstenholme primes}}

Sloane, N.J.A. Sequences A001008, A007406, A007408, A088164, 
A034602, A177783, in "The On-Line Encyclopedia of Integer Sequences." 
(published electronically at www.research.att.com/~njas/sequences/).

\vspace{3mm}

\centerline{\bf{C) List of papers/authors ordered by years of publications}}

\vspace{2mm}

{\small
\cite[1801]{ga} C.F. Gauss

\cite[1819]{bab} C. Babbage

\cite[1851]{k} E.E. Kummer

\cite[1862]{w} J. Wolstenholme

\cite[1877-1878]{l}   E. Lucas

\cite[1889]{le} C. Leudesdorf

\cite[1900]{gl1} J.W.L. Glaisher

\cite[1900]{gl2} J.W.L. Glaisher

\cite[1901]{gl3} J.W.L. Glaisher

\cite[1930]{sv} E.T. Stafford and H.S. Vandiver

\cite[1933]{ch1} S. Chowla

\cite[1934]{ch} S. Chowla

\cite[1934]{hw2} G.H. Hardy and  E.M. Wright

\cite[1935]{gup} H. Gupta

\cite[1937]{ra} R. Rao

\cite[1938]{leh} E. Lehmer

\cite[1941]{se} D. Segal/H.W. Brinkman 

\cite[1952]{bs}  V. Brun, J.O. Stubban, J.E. Fjeldstad, R. Tambs Lyche,
K.E. Aubert, W. Ljunggren and E. Jacobsthal

\cite[1954]{ca1} L. Carlitz

\cite[1955]{dp} H.J.A. Duparc and W. Peremans

\cite[1968]{kaz} G.S. Kazandzidis

\cite[1969]{kaz2} G.S. Kazandzidis

\cite[1974]{tr} Yu.A. Trakhtman

\cite[1975]{jo} W. Johnson

\cite[1979]{ap} R. Ap\'ery  

\cite[1982]{ge} I. Gessel

\cite[1983]{ge2} I. Gessel

\cite[1988]{bau} F.L. Bauer

\cite[1988]{gar} A. Gardiner

\cite[1988]{ro} N. Robbins

\cite[1989]{zhan} S. Zhang

\cite[1990]{ba} D.F. Bailey

\cite[1991]{el} A. Eswarathasan and E. Levine

\cite[1992]{ta} Y. Tan

\cite[1993]{bcem} J. Buhler, R. Crandall, R. Ernvall and T. Mets\"{a}nkyl\"{a}

\cite[1994]{a} E. Alkan

\cite[1994]{bo} D.W.  Boyd

\cite[1995]{cla} W.E. Clark

\cite[1995]{che} K.Y. Chen

\cite[1995]{hon} S. Hong

\cite[1995]{m} R.J. McIntosh

\cite[1995]{rz} A. Robert and M. Zuber

\cite[1997]{b} M. Bayat  

\cite[1997]{gr} A. Granville

\cite[1998]{bh} D. Berend and J.E. Harmse

\cite[1998]{ge1} I. Gessel

\cite[1998]{wa} L. Washington

\cite[1999]{an} G.E. Andrews

\cite[1999]{kw} W. Kimball and W. Webb

\cite[1999]{sl1}  I.Sh. Slavutskii

\cite[1999]{sh1} Q. Sun and S. Hong

\cite[2000]{hon2} S. Hong

\cite[2000]{sl2}  I.Sh. Slavutskii

\cite[2000]{s1} Z.H. Sun

\cite[2001]{sh2} Q. Sun and S. Hong

\cite[2001]{tw} V. Trevisan and K.E. Weber

\cite[2002]{cg}  T.X. Cai and A. Granville

\cite[2002]{sl3}  I.Sh. Slavutskii

\cite[2004]{gru} D.B. Gr\"{u}nberg

\cite[2006]{chd1} M. Chamberland and K. Dilcher

\cite[2006]{zy} H. Zheng and S. Yang

\cite[2007]{hon} S. Hong

\cite[2007]{mr} R.J. McIntosh and E.L. Roettger

\cite[2007]{pan} H. Pan

\cite[2007]{sp} L.L. Shi and H. Pan

\cite[2007]{sd} Z.W. Sun and D.M. Davis

\cite[2007]{sw} Z.W.  Sun and D. Wan

\cite[2007]{z2} J. Zhao

\cite[2007]{zc} X. Zhou and T. Cai

\cite[2008]{dil} K. Dilcher  

\cite[2008]{ht} C. Helou and G. Terjanian

\cite[2008]{s2} Z.H. Sun

\cite[2008]{pan2} H. Pan

\cite[2008]{z1} J. Zhao

\cite[2009]{chd2} M. Chamberland and K. Dilcher

\cite[2009]{pan3}  H. Pan

\cite[2010]{bls} K.A. Broughan, F. Luca and I.E. Shparlinski

\cite[2010]{t} R. Tauraso

\cite[2010]{ylz} J. Yang, Z. Li and F. Zhang

\cite[2010]{z3} J.  Zhao

\cite[2011]{me1} R. Me\v{s}trovi\'c

\cite[2011]{me2} R. Me\v{s}trovi\'c

\cite[2011]{me3} R. Me\v{s}trovi\'c

\cite[2011]{me4} R. Me\v{s}trovi\'c

\cite[2011]{me5} R. Me\v{s}trovi\'c

\cite[2011]{str} A. Straub

\cite[2011]{syl} Y. Su, J. Yang and S. Li

\cite[2011]{su1} Z.W.  Sun

\cite[2011]{su3} Z.W.  Sun

\cite[2011]{su2} Z.W.  Sun

\cite[2011]{t2} R. Tauraso

\cite[2011]{ya} J. Yang
      \end{document}